\documentclass[11pt,leqno]{amsart}
\usepackage{amscd}
\usepackage{amssymb}
\usepackage{amsfonts}
\usepackage{latexsym}
\usepackage{verbatim}
\usepackage{amsthm}
\usepackage{mathrsfs}
\usepackage[small,nohug,heads=littlevee]{diagrams}
\diagramstyle[labelstyle=\scriptstyle]

\usepackage{url,euscript,graphics,graphicx,a4wide}

\usepackage[british]{babel}
\usepackage[T1]{fontenc}

\theoremstyle{remark}     
\newtheorem{remark}{Remark}[section]\newtheorem*{remark*}{Remark}

\def\sideremark#1{\ifvmode\leavevmode\fi\vadjust{\vbox to0pt{\vss
\hbox to 0pt{\hskip\hsize\hskip1em%
\vbox{\hsize2cm\tiny\raggedright\pretolerance10000%
\noindent #1\hfill}\hss}\vbox to8pt{\vfil}\vss}}}%

\newcommand{\be}{\begin{equation}}
\newcommand{\ee}{\end{equation}}

\def\GG{\mathscr{G}}

\def \T{\mathrm{T}}

\newcommand{\rP}{\mathrm P}
\newcommand{\rQ}{\mathrm Q}
\newcommand{\rR}{\mathrm R}

\newcommand{\iso}{\cong}

\newcommand{\lto}{\longrightarrow}

\newcommand{\hook}{\lrcorner\,}

\renewcommand{\geq}{\geqslant}\renewcommand{\leq}{\leqslant}
\renewcommand{\H}{\mathbb{H}}

\newcommand{\La}{\Lambda}
\newcommand{\hodge}{{*}}

\DeclareMathOperator{\End}{\textsl{End}\,}

\newcommand{\ol}{\overline}
\newcommand{\q}{\quad}
\newcommand{\ba}{\begin{array}}\newcommand{\ea}{\end{array}}
\newcommand{\vs}{\vphantom{\int_W^M}}

\numberwithin{equation}{section}
\theoremstyle{plain}
\newtheorem{theorem}{Theorem}[section]
\newtheorem{definition}[theorem]{Definition}
\newtheorem{lemma}[theorem]{Lemma}
\newtheorem{prop}[theorem]{Proposition}
\newtheorem{cor}[theorem]{Corollary}

\newtheorem{ex}[theorem]{Example}
\newtheorem{exs}[theorem]{Examples}

\newcommand{\SO}{\mathrm{SO}}
\newcommand{\GL}{\mathrm{GL}}
\newcommand{\so}{\mathfrak{so}}
\newcommand{\U}{\mathrm{U}}
\newcommand{\SU}{\mathrm{SU}}
\newcommand{\su}{\mathfrak{su}}
\newcommand{\psu}{\mathfrak{psu}}
\newcommand{\Sp}{\mathrm{Sp}}
\newcommand{\ssp}{\mathfrak{sp}}
\newcommand{\G}{\mathrm{G}}
\newcommand{\g}{\mathfrak{g}}
\newcommand{\sW}{\mathscr{W}}
\newcommand{\sG}{\mathscr{G}}
\newcommand{\sg}{\mathrm{\bf g}}

\newcommand{\Spin}{\mathrm{Spin}}
\newcommand{\rS}{\mathrm S }
\newcommand{\rE}{\mathrm {E}}
\newcommand{\rA}{\mathrm A}

\newcommand{\aaa}{\alpha}
\newcommand{\bbb}{\beta}
\newcommand{\R}{\mathbb{R}}
\newcommand{\C}{\mathbb{C}}
\newcommand{\Z}{\mathbb{Z}}

\newcommand{\Span}{\mathrm{Span}}
\newcommand{\rH}{\mathrm{H}}
\newcommand{\PSU}{\mathrm{PSU}(3)}

\newcommand{\M}{\mathrm{M}}
\begin{document}
\title[$\SO(3)$-structures on $8$-manifolds]
{$\boldsymbol{\SO(3)}$-structures on $\boldsymbol{8}$-manifolds}
\author[S.G.\ Chiossi,\; \'O.\ Maci\'a]{Simon G. Chiossi$^{1)}$,\; \'Oscar Maci\'a$^{1,2)}$}

\address{\noindent
$1)$ Dipartimento di Matematica
\\ Politecnico di Torino \\ Corso Duca degli Abruzzi 24, 10129 Torino,
Italy}
\email{oscarmacia@calvino.polito.it}
\email{simon.chiossi@polito.it}
\address{\noindent
$2)$ Departamento de Geometr\'ia y Topolog\'ia\\
Universidad de Valencia\\
C. Dr. Moliner, S/N, 96100 Burjassot, Valencia, Spain}
\email{oscar.macia@uv.es}

\thanks{This work was supported in part by the project {\sc Prin} \oldstylenums{2007} ``Differential Geometry and Global Analysis'' of {\sc Miur} (Italy), and
by the research grant {\sc VALi+d APOSTD/2010/044} from the Generalitat Valenciana (Spain).}
\subjclass[2010]{primary 53C10; secondary 53C15, 53C26}
\keywords{$\SO(3)$-structure, Riemannian $8$-manifold, intrinsic torsion,
quaternion-Hermitian, $\PSU$}

\begin{abstract}
We study Riemannian $8$-manifolds with an infinitesimal action of $\SO(3)$
by which each tangent space breaks into irreducible spaces of dimensions $3$ and $5$.
The relationship with quaternionic, almost product- and $\PSU$-geometry is thoroughly explained using representation-theoretical arguments.
\end{abstract}

\maketitle
\frenchspacing

\tableofcontents

\section{Introduction}

\noindent
The purpose of this article is to study geometric structures on
eight-dimensional Riemannian manifolds admitting a special infinitesimal
action of the Lie group $\SO(3).$ Our study represents the first flicker of the unfolding of a larger theory, as we shall have time to argue; it follows a long trail initiated
by A.Gray and L.Hervella \cite{Gray}, who recognised the intrinsic torsion
 as the killer app to interpret, and handle, Riemannian manifolds $\{\M^{n},g\}$ whose principal bundle of oriented orthonormal frames reduces to a bundle with structure group $\G\subset \SO(n)$, aka $\G$-structures.

Among the flurry of $\G$-structures considered ever since -- for a comprehensive survey of which the adoption of \cite{Agricola:Srni} is recommended -- we are concerned with those for which $\G=\SO(3)$.
This Lie group can act reducibly on the tangent spaces of a manifold, a situation sorted out years ago by  E.Thomas \cite{Thomas:vf} and M.Atiyah \cite{Atiyah:vf}, or irreducibly, in which case the general stance can be found in \cite{Friedrich}.
 The first concrete results about an $\SO(3)$-irreducible action were obtained in \cite{Nurowski} and \cite{FinoChiossi} on manifolds of dimension $5$. The former represented a true breakthrough, for M.Bobi\'enski and P.Nurowski proved that the r�le of the intrinsic torsion of an irreducible $\SO(3)$-manifold is played
  by a tensor defined naturally by the representations.
 Further advancement was made by I.Agricola et al. \cite{Friedrich2}, in relationship to the topological obstructions to existence in dimension five in particular.\\

As a matter of fact $\SO(3)$  is a quotient of the
Lie group $\SU(2)$, and as such plays a role, in real dimensions $4k$, analogous to that of the circle group $\U(1)$ in complex geometry.
This paper aims at describing $\SO(3)$-structures on Riemannian manifolds
$\{\M^{8},g\}$ for which each tangent space decomposes in two summands
$$
\T_{p}\M^{8} = V\oplus W,
$$
of dimensions three and five, each of which is irreducible under $\SO(3)$.
This is achieved by fixing a representation  of $\SO(3)$ on $\R^{8}$,  whose image we call $\GG$.
The reason for concentrating on dimension eight and focusing on that particular action is that
$$
\GG=\SO(3)
$$
is contained in other intermediate Lie subgroups $\G$ of $\SO(8)$ of a certain interest, namely
$$
\SO(3)\times \SO(5),\q  \PSU,\q \Sp(2)\Sp(1).
$$
Each group of this triad is known to give rise to exciting geometric properties, see \cite{Naveira}, \cite{Hitchin:stable, Witt:triality} and \cite{Salamon:dgH,Swann} respectively, and is remarkably
obtained for free when imposing the $\GG$-action. What kickstarted the present work was
the realisation that there is a deep relationship between $\sG$ and almost quaternion-Hermitian geometry, as noticed also by A.Gambioli \cite{Gambioli:latent}. The article \cite{Macia}, to which this is a sequel of sorts, made this
idea concrete by studying the quaternionic geometry of $\SU(3)$, focusing in particular on the class of so-called nearly quaternionic structures, the closest quaternionic kin to non-complex nearly K\"ahler manifolds. Whereas the customary approach would consider a Lie group's bi-invariant metric properties, in that case there is one action of $\SU(3)$ on the left, and another action of $\SO(3)$ on the right.

One important point is to show explicitly how the `subordinate' $\G$-structures relate to one another, and especially how they determine a $\GG$-action. To this end we prove in Theorem \ref{main} that not only $\GG$ is the triple intersection of the aforementioned groups, as expected, but one can recover it by using just two of the three. This result is rephrased with more algebraic flavour by Theorem  \ref{gig1}, according to which, if we let $\g\subset \so(8)$ be the Lie algebra of $\sG$, then
 $$
\g^\perp = \frac{\so(3)\oplus\so(5)}{\g} +
 \frac{\psu(3)}{\g} +
\frac{\ssp(2) \oplus \ssp(1)}{\g},
$$
and the sum of any two quotients is isomorphic to the orthogonal complement of the third algebra.

Triality, a distinctive feature of dimension $8$, manifests itself spectacularly by endowing the space $\SO(8)/\PSU$, parametrising reductions of the structure group of $\T\M^{8}$ to $\PSU$, with the structure of a $3$-symmetric space (Remark \ref{remtri1}).

After a brief diversion on topology (section \ref{sec:top}), which serves the purpose of finding obstructions to the existence of $\GG$-actions, we dwell into the theory of $\G$-structures.
We determine the intrinsic-torsion space for $\GG$ by decomposing it into irreducible $\GG$-modules  (Proposition \ref{SO(3)it}); the core observation is that, of its $200$ dimensions,   exactly $3$ arise from $\sG$-invariant tensors that we can describe explicitly. We then compare this  decomposition with the spaces of $\G$-intrinsic torsion relative to the common subgroup $\GG$, culminating in Proposition \ref{prop:torsion}.

A further refinement is presented in the last section \ref{sec:invariant}.
There are six exterior forms invariant under $\GG$, two $3$-forms, two $4$-forms and two $5$-forms,  defining the reduction.
If one restricts to the case in which the intrinsic torsion is $\GG$-invariant, Theorem \ref{prop:cases} guarantees these three pairs will fit into an invariant deRham complex
$$
(\La^{3}{\mathrm{T^{*}M}})^{\GG} \stackrel{d}{\lto}
 (\La^{4}{\mathrm{T^{*}M}})^{\GG} \stackrel{d}{\lto}
(\La^{5}{\mathrm{T^{*}M}})^{\GG}.
$$
that governs the geometry entirely.
\\

We should emphasize that this research area is very much in its infancy. We have concentrated on the Lie-theoretical point of view embodied by the intrinsic torsion, and put
little stress on other aspects that should deserve a separate treatment, such as  higher-dimensional Riemannian $\SO(3)$-manifolds.\\

\noindent
{\it The authors are indebted to Simon Salamon for sharing his insight,
and wish to thank Thomas Friedrich and Andrew Swann for enlightening
conversations.
}
%


\section{$\SO(3)$ and subordinate structures}
\noindent
We recall the salient points of $\Sp(1)$-representations, since the chief technique for decomposing exterior forms employs the weights of the action of $\ssp(1)\otimes\C=\mathfrak{sl}(2,\C)$. Let $\rH$ be the basic representation of $\Sp(1)$ on $\C^{2}=\mathbb{H}$ given by left-multiplying column vectors by matrices. The $n$th symmetric power of $\rH$ is an irreducible representation of $\C^{2n}=\mathbb{H}^{n}$ written
$$
\rS^{n}=\rS^{n}(\rH).
$$
Every irreducible $\Sp(1)$-module is of this form for some nonnegative $n\in\Z$, and is the eigenspace of the Casimir operator with eigenvalue $-n(n+2)$. It can be decomposed into weight spaces of dimension one under the action of the Cartan subalgebra $\mathfrak{sl}(2,\C)$, making $\rS^{n}$ the unique irreducible representation with highest weight $n$. Since the weight is a homomorphism from tensor products of representations to the additive integers, weight-space decompositions can be used to determine tensor-, symmetric and skew-symmetric products of modules, as in the Clebsch-Gordan formula
$$
\rS^{n}\otimes\rS^{m}\iso\rS^{n+m}\oplus\rS^{n+m-2}\oplus\ldots\oplus
\rS^{n-m+2}\oplus\rS^{n-m}
$$
when $n\geq m$.

The Lie group $\SO(3)$ is the $\Z_{2}$-quotient of the simply connected covering $\Sp(1)$, and as such its complex representations coincide with the aforementioned ones. Thus by identifying
$\SO(3)$-modules with the above $\rS^{\lambda}$  endowed with a real structure, we will allow  ourselves to switch tacitly from complex to real throughout this work. For example $V\cong \R^3$,  the fundamental representation of $\SO(3)$, will be considered the same as $\rS^{2}\rH$, and the traceless symmetric product $\rS^{2}_{0}V\iso \rS^{4}\rH$ will have dimension five.\\

A natural way to manufacture an $\SO(3)$-structure on a Riemannian manifold
$\{\M^{8},g\}$ is to fix a  homomorphism
\begin{equation}\label{q-inclusion}
\SO(3)\lto \SO(8),
\end{equation}
whose image will be denoted $\GG$,
that breaks the tangent space at each point $p$ of $\M$ into three- and five-dimensional summands
\begin{equation}\label{3+5=8}
\T_{p}\M^{8} = V\oplus \rS^{2}_{0}V=V\oplus W.
\end{equation}
Note that $W$ is in fact the irreducible $\SO(3)$-module used in \cite{Nurowski, FinoChiossi, Friedrich} to study $\SO(3)$-manifolds of dimension $5$ modelled on the Riemannian symmetric space $\SU(3)/\SO(3)$. The story there went as follows: there is a subgroup inside $\SO(5)$ isomorphic to $\SO(3)$ acting in an irreducible fashion, by which $\su(3)=\so(3)+iW$, and $W\iso\R^{5}$ is identified with the space of symmetric and trace-free $3\times 3$ matrices $\rS^{2}_{0}\R^{3}$.
Although the $\GG$-action prescribed by \eqref{3+5=8} is certainly not irreducible, it is the unique action of $\SO(3)$ that decomposes $\R^{8}$ into non-trivial $\GG$-irreducible subspaces, and as such it is worth studying.

A more interesting way to view the splitting is to rewrite it as
\begin{equation}\label{splitting}
\T_{p}\M^{8}\otimes_{\R}\C\cong \rS^2\rH\oplus\rS^4\rH.
\end{equation}
If a Riemannian $8$-manifold admits such an infinitesimal action of $\SO(3)$, the endomorphisms of the complexified tangent bundle $\T\M^{*}_{c}\otimes\T\M_{c}$ include a factor $\rS^2$, which is precisely the Lie algebra $\g\iso \so(3)$ of $\sG$.

There is a second reason for concentrating on dimension eight, and focusing on the particular action described by \eqref{splitting}. The embedding of $\SO(3)$ actually factors through other Lie groups of interest, refining \eqref{q-inclusion} to this diagramme:
\begin{diagram} \label{triangle}
\GG & & \rTo & & \SO(8) \\
& \rdTo & & \ruTo & \\
 & & \G & &  \\
\end{diagram}
where $\G$ will be one of the following  subgroups of $\SO(8)$:
\be\label{eq:3G}
\Sp(1) \Sp(2)=\Sp(1) \times_{\Z_{2}}\Sp(2),\q \PSU={\SU(3)}/\Z_{3}, \q \SO(3)\times \SO(5).
\ee
This enables us to understand the mutual relationship between an $\SO(3)$-structure and any of these $\G$-structures. The main point to this section is to prove
\begin{prop}\label{factors}
A $\sG$-structure \eqref{splitting} on a  Riemannian eight-manifold $\{\M^8,g\}$ induces altogether an almost product structure, a $\PSU$-structure and an almost quaternion-Hermitian structure.
\end{prop}
\proof
The argument will be broken up in cases corresponding to the above Lie groups $\G$.
\begin{itemize}
\item[(i)] $\G=\Sp(2)\Sp(1)$.
What we will actually prove here is the intermediate `diagonal-type' inclusion
$$\sG\subset \Sp(1)_+\times\Sp(1)_-\subset \Sp(2)\Sp(1),$$
whereby the first factor is the irreducible $\Sp(1)$ inside $\Sp(2)$ (recall the latter is the universal cover of $\SO(5)$), while the second is just the identity map's image.
Similar $\pm$ labelling will be used to identify the representations.

At each point of $\M$ the complexified tangent space $\T_c=\rS^{2}\oplus\rS^{4}$
is quaternionic, since the action of $\Sp(1)$ on the fundamental $\Sp(2)$-representation
$\rE\cong \C^4$  with highest weight $(1,0)$ ensures $\rE=\rS^3\rH$. This means  that
\eqref{splitting} is isomorphic to
$$
\T_c\cong \rS^3_{+}\rH\otimes \rH=\rE\otimes\rH,
$$
if one takes, by convention, $\rH=\rS^{1}_{-}\rH$.
The orthogonal Lie algebra decomposes under $\Sp(1)_+\times \Sp(1)_-$ as
\be\label{so8aqh}
\ba{rccccc}
\so(8,\C)= \Lambda^2(\rS^3_+ \rS^1_-)= &
      \left(\rS^6_+ \oplus  \rS^2_+ \right) & \oplus & \left(\Lambda^2_0(\rS^3_+)\otimes \rS^2( \rS^1_-)\right) & \oplus & \rS^2_- \\[1mm]
=& \ssp(2)_{+} &\oplus & (\ssp(2)\oplus\ssp(1))^\perp &\oplus &  \ssp(1)_-.
\ea
\ee
The group $\sG$ then sits inside $\Sp(1)_+\times\Sp(1)_-$ and its
 Lie algebra $\g$ corresponds to an $\rS^2$-module in
 $\rS^2_+\oplus\rS^2_-=\ssp(1)_+\oplus\ssp(1)_-.$
 \smallbreak


\item[(ii)] $\G=\PSU$.
The corresponding Lie algebra $\psu(3)=\su(3)$ splits naturally as $\rS^2\oplus\rS^4$ if we think of
a $3\times 3$ Hermitian matrix as the Cartan sum of a skew-symmetric matrix (in $\rS^2=\so(3)$)
and a purely imaginary symmetric matrix with no trace (whence $W=\rS^{4}$).
\smallbreak


\item[(iii)] $\G=\SO(3)\times \SO(5)$.
The argument is formally the same of case (i), by virtue of the universal covering
$\Sp(1)\times\Sp(2)\stackrel{4:1}{\lto}\SO(3)\times \SO(5)$.
In two words, the block-diagonal embedding $\SO(3)\times \SO(5)\subset \SO(8)$
 decomposes $\so(8)\cong \Lambda^2(W\oplus V)
= \so(5)_W\oplus \left(W\otimes V \right)\oplus \so(3)_V$, reflected in \eqref{so8aqh}.
\hfill qed
\end{itemize}
\begin{ex}\label{ex:oscar}
A class of $\sG$-structures on $\M^{8}=\SU(3)$, giving a non-integrable analogue of quaternion-K\"ahler geometry, was studied in \cite{Macia}. The action of $\g$ on $\su(3)$ described there
  spawned an almost quaternion-Hermitian structure of class $\sW_1^{AQH}\subset \sW_{1+4}^{AQH}$ in the Cabrera-Swann terminology; the fact that the local K\"ahler forms generated a differential ideal with coefficients sitting in a traceless, symmetric matrix of $1$-forms, was the main reason to start the study of $\sG$-manifolds.
\end{ex}

\noindent
Before passing to the next section, a few credits to complete the overall picture.
\smallbreak

(i)
Recall $\Sp(2)\Sp(1)$ is one of Berger's holonomy groups \cite{Berger}, see
\cite[ch. 5]{Besse2} for an account of the theory.
To S.Salamon \cite{Salamon:QK} we owe the `EH' formalism used throughout the paper, while
for the exhaustive description of the intrinsic torsion of  $\Sp(2)\Sp(1)$-geometry one should refer to
\cite{Swann:symplectique,Cabrera2007}.
\smallbreak

(ii)
$\PSU$-structures were borne in on N.Hitchin's programme on special geometry
\cite{Hitchin:stable} and were thoroughly explored by F.Witt \cite{Witt:triality}.
The latter and Puhle's article \cite{Puhle} give an accurate description of the intrinsic
torsion of $\PSU$-structures, and our results are meant to complement those and enable to
gain solid, intuitive understanding of the matter.
\smallbreak

(iii)
The classification of almost product structures begun with the work of
A.Naveira \cite{Naveira}; deep geometric consequences and many examples were discussed
by O.Gil-Medrano \cite{GM} and V.Miquel \cite{Miquel}.


\section{Intersection theorem}
\label{intersection}

\noindent
Having shown that $\SO(3)$ is a subgroup of all three of
$\SO(3)\times\SO(5)$,\;$\Sp(2)\Sp(1)$,\;$\PSU$, we know in particular it
 lies in their intersection. The threefold aspect of the theory is pivotal for understanding
 this and many other features, and the aim of this section is to prove that there is some
 redundancy: in other words, any two are enough to guarantee the retrieval of $\GG=\SO(3)$.
 More formally,
\begin{theorem}\label{main}
Let $\sG\iso\SO(3)$ be the subgroup of $\SO(8)$ acting infinitesimally on
a Riemannian $8$-manifold $\{\M,g\}$ by decomposing tangent spaces like
$$
\T_{p}\M= V\oplus \rS^2_0V,
$$
 where $V\cong\R^3$ is the fundamental representation. Then
\begin{enumerate}
\item $\sG=\bigl(\SO(3)\times\SO(5)\bigr)\cap\PSU$,\smallbreak
\item $\sG=\PSU\cap\Sp(2)\Sp(1)$,\smallbreak
\item $\sG=\Sp(2)\Sp(1)\cap\bigl(\SO(3)\times\SO(5)\bigr).$
\end{enumerate}
\end{theorem}

\begin{figure}[h!]
\centering
  \includegraphics[width=6cm]{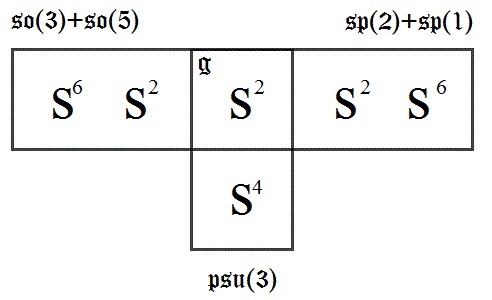}
\caption{$\so(8)$ as the sum of $\so(3)\oplus\so(5),\;\psu(3),\;\ssp(2)\oplus\ssp(1),$ neatly intersecting in $\g.$}\label{fig1}
\end{figure}

The proof, sketched in fig. \ref{fig1}, consists of the ensuing series of lemmas and is minutely carried out at the level of the corresponding Lie algebras.

\begin{lemma}
Retaining the notation of \eqref{3+5=8},
$$
\g = \bigl(\so(3)_V\oplus\so(5)_W\bigr)\cap \su(3).
$$
\end{lemma}

\proof The  right-hand-side intersection, say $\sg$, contains the diagonal $\g.$
Now, $\g$ is sitting diagonally in $\so(3)_V\oplus\so(5)_W$, so
the combination  $\g\oplus\so(3)_V$ would generate $\so(5)_W$
if $\su(3)$ were to contain $\g$ together with $\so(3)_V$. Since
 $\so(5)_W$ is not a subalgebra of $\su(3)$, we have $\g=\sg.$ \hfill qed

\begin{lemma}\label{psu&quat}
 We have
$$
\g = \su(3)\cap \bigl(\ssp(2) \oplus \ssp(1)\bigr).
$$
\end{lemma}

\proof
Let $\sg$ be the intersection on the right, and recall  $\g\subseteq \sg$. Moreover,
$\g$ sits diagonally in $\ssp(1)_+\oplus \ssp(1)_-\subset\ssp(2)\oplus \ssp(1)$.
The fundamental representation  $\C^3$ of $\su(3)$ corresponds either to $\C^3=\rS^2$ or $\C^3=\rH\oplus\C$, if viewed as a $\g$-module. But since $\su(3)$ is the space of traceless endomorphisms of $\C^3$, the two instances regard $\End_0(\rS^2)=\rS^4\oplus\rS^2$ or $\End_0(\rH\oplus\C)=\rS^2\oplus 2\rH\oplus \C$ respectively. Neither contains a $3+3$-dimensional submodule that can be identified with $\ssp(1)_+\oplus \ssp(1)_-$, so $\ssp(1)_+\oplus \ssp(1)_-\nsubseteq \su(3)$ and the $\rS^2$ factor of the first instance must coincide with $\g.$ \hfill qed\\

\begin{lemma}
In the notation of \eqref{3+5=8},
$$
\g=\bigl(\ssp(2) \oplus \ssp(1)\bigr)\cap \bigl(\so(3)_V \oplus \so(5)_W\bigr).
$$
\end{lemma}
\proof

\indent Let $\sg$ denote the intersection, as usual. As we know that $\g\subset\sg$, let us assume by contradiction $\g\neq \sg$.

The following possibilities arise when looking at $\sg\subset \ssp(2)\oplus\ssp(1)$: either $(i)$ $\ssp(1)\subseteq \sg,$ or $(ii)$ $\ssp(1)\nsubseteq \sg.$

Case $(i)$: $\sg=\ssp(1)\oplus \mathfrak{k},$ with $\mathfrak{k}=\sg\cap\ssp(2).$ Let $\G$ and $\mathrm K$ denote the connected Lie groups of $\sg,\;\mathfrak{k}$.  Complex tangent spaces decompose as $\T_c = \rE\rH$, and as $\ssp(1)\subseteq \sg,$ the $\Sp(1)$-irreducible representation $\rH$ is irreducible also under $\G.$
  Similarly $\rE$ can be seen as a $\mathrm K$-module, since $\mathfrak{k}\subseteq \ssp(2)$; in contrast to $\rH$, however, it can be either reducible or irreducible, depending on whether  $\mathfrak k\subset \ssp(1)_+\subset \ssp(2)$ (case $(i.i)$) or  $\ssp(1)_+\subseteq \mathfrak k\subseteq \ssp(2)$ (case $(i.ii)$).

$(i.i)$: write $\rE=\oplus_i \rA_i$ as a sum of irreducible $\mathrm K$-modules, so that $\T_c =\left(\oplus_i\rA_i\right)\otimes \rH=\oplus_i\left(\rA_i\otimes \rH\right)$ is a sum of $\G$-irreducible terms of even dimension. The latter fact clashes with the dimensions of $\R^5\oplus\R^3$ arising from
$\sg\subset \so(5)\oplus\so(3)$.

$(i.ii)$: $\rE$ is $\mathrm K$-irreducible, making $ \T_c = \rE\rH$ irreducible under $\G$.
But again, \eqref{splitting} contradicts irreducibility.

 Case $(ii)$: as $\ssp(1)\nsubseteq \sg$, we write $\sg\subseteq \mathfrak k \oplus\ssp(1),$ where now $\mathfrak k$ is the projection of $\sg$ to $\ssp(2).$ But since $\g$ acts diagonally, $\ssp(1)_+\subseteq\mathfrak k \subseteq \ssp(2),$ and $\mathfrak k$ is a $\Sp(1)_+$-representation inside $\ssp(2).$
Considering the decomposition of $\ssp(2)=\rS^6_+\oplus \rS^2_+$ in irreducible $\Sp(1)_+$-modules, with $\rS^2_+\cong\ssp(1)_+$, we face another dichotomy: either $(ii.i)$ $\mathfrak k = \ssp(1)_+$, implying that $\g\oplus\ssp(1)_+\subseteq \sg,$ or $(ii.ii)$ $\mathfrak k = \rS^6_+$ and $\g\oplus \rS^6_+\subseteq\sg.$

$(ii.i)$: as $\g$ sits diagonally in $\ssp(1)_+ \oplus \ssp(1)_-,$ the subalgebra $\g \oplus \ssp(1)_+ \subseteq \sg$ would detect an $\ssp(1)_-=\ssp(1)$ inside $\sg,$ against the general hypothesis $(ii).$

$(ii.ii)$: $\rS^6_+$ is not a subalgebra of $\ssp(2),$ so the Lie bracket of $\g \oplus\rS^6_+\subseteq \sg $ would produce an $\ssp(2)$-term inside $\sg.$ By the same argument as before, the diagonal $\g\subset\ssp(2)\oplus\ssp(1)$ would make $\g\oplus\ssp(2)$ generate an $\ssp(1)$ in $\sg,$ again against  $(ii).$

 Both $(i),\;(ii)$ disproving the initial assumption, we conclude that $\g=\sg.$\hfill qed\\

We can rephrase Theorem \ref{main} in a perhaps-more-eloquent fashion
\begin{theorem}\label{gig1}
Let $\g_i,\;i=1,2,3$, denote the Lie algebras of the groups $\SO(3)\!\times\!\SO(5)$, $\PSU$, $\Sp(2)\Sp(1)$,
$\g_i^\perp$ the complements in $\so(8)$ and $\g$ the Lie algebra of $\sG=\SO(3).$ Then
\begin{eqnarray*}
\g_i^\perp & = &  (\g_j/\g) \oplus (\g_k/\g),\qquad i\neq j \neq k=1,2,3\\
\g^\perp  & = & \bigoplus_{i=1}^3 (\g_i/\g).
\end{eqnarray*}
\end{theorem}
\proof
As $\so(8)=\g_i \oplus \g_i^\perp$ and $\g_i=\g \oplus (\g_i/\g)$, the assertion is
a straightforward consequence of a dimension count plus
$\g_i\cap\g_j=\g$ (Theorem \ref{main}). \hfill qed\\

For the sake of clarity, and for later use, here are the explicit modules involved.
Start from the irreducible $\sG$-decomposition
\be\label{eq:2forms}
\so(8) \cong \Lambda^2(\rS^2\oplus\rS^4) = 2\rS^6\oplus
  \rS^4\oplus 3\rS^2 \cong \g \oplus \left(2\rS^6\oplus\rS^4\oplus 2\rS^2\right),
\ee
where $\g\cong\rS^2,\;\g^\perp=(2\rS^6\oplus\rS^4\oplus 2\rS^2).$
Read in `VW' terms, that tells
 $(\so(3)_V\oplus\so(5)_W)^\perp\cong\rS^2_V\otimes\rS^4_W$ as
 $\SO(3)_V\times \SO(3)_W$-modules. Reducing to $\sG$-modules by taking
 the diagonal embedding (forgetting where the terms come from, ie dropping the subscripts)
 leads to
$$
\bigl(\so(3)\oplus\so(5)\bigr)^\perp \cong \rS^6 \oplus \rS^4\oplus \rS^2
$$
hence $\bigl(\so(3)\oplus\so(5)\bigr)/\g = \rS^2\oplus\rS^6.$

Similarly, starting from (\ref{so8aqh}) and
 identifying $\rS^{\lambda}_+\cong\rS^{\lambda}_-$, we obtain
$$
\bigl(\ssp(2)\oplus\ssp(1)\bigr)^\perp \cong
  \bigl(\Lambda^2_0(\rS^3)\otimes\rS^2(\rS^1)\bigr)\cong
  \rS^4\otimes\rS^2 = \rS^6\oplus\rS^4 \oplus \rS^2
$$
with $\bigl(\ssp(2)\oplus\ssp(1)\bigr)/\g = \rS^6 + \rS^2.$

Finally, it is well documented \cite{Macia} that the adjoint
representation of $\su(3)$ in $\so(8)$, decomposed under
$\sG$, coincides with $\rS^2\oplus\rS^4.$ Hence
\be\label{psu-perp}
\su(3)^\perp = 2\rS^6 \oplus 2\rS^2
\ee
and $\su(3)/\g \cong \rS^4.$

 \begin{figure}[h!]
 \centering
 \includegraphics[width=6cm]{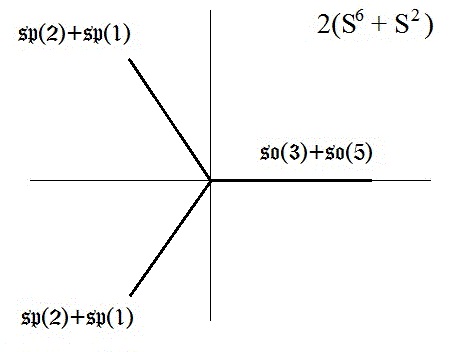}
 \caption{The space $2(\rS^2\oplus\rS^6)$ rotates by $2\pi/3$ incarnating at each turn a pair of   subalgebras (bold lines), while $\psu(3)$ remains fixed (centre of symmetry).}\label{fig-triality}
 \end{figure}

\begin{remark} [\cite{Witt:triality}]\label{remtri1}
The triality automorphism of $\R^{8}$,  in some loose sense, plays a similar role
to that of self-duality in dimension $4$. Its stabiliser is $\G_2\subset \Spin(8)$,
and although the exceptional Lie group does not contain $\sG$, but rather a diagonal
$\SO(3)$ embedded in $\SO(4)\subset\G_2$,  the decomposition of $\so(8,\C)$ detects
the existence of an automorphism of order three (fig. \ref{fig-triality}) that
permutes the three $\rS^2$ by rotating the two copies of $\rS^6\oplus\rS^{2}$ by $2\pi/3$. This automorphism  has $\PSU\times\Z_2$ as fixed-point set in $\Spin(8)$
(and $\PSU$ in $\SO(8)$).

The $20$-dimensional manifold $\SO(8)/\PSU$ is a $3$-symmetric space \cite{GrayWolf}
whose tangent space at each point is isomorphic to $\rS^2\C^3\oplus\rS^2\ol{\C^3}.$ This is given by \eqref{psu-perp} and is real-irreducible under $\PSU$, whereas over $\SO(3)$ it can be decomposed as a sum $(\rS^2\oplus\rS^6)\oplus f_*(\rS^2\oplus\rS^6),$ where $f$ is the induced isometry of order three.
 If the  sum of the first two terms corresponds to $\so(3)\oplus\so(5),$  then  $f_*(\rS^2\oplus\rS^6)$ is $\ssp(2)\oplus\ssp(1)$. The third space required to `visualise'  triality is the image of the square  $f_*\circ f_*$, and is only apparently missing. It is another quaternionic subalgebra fixed by an `anti-selfdual' four-form \cite{Macia}, and we will return to it later (see \eqref{Omega'}).
\end{remark}


\section{Topological observations}\label{sec:top}
\noindent
This section marks a slight detour intended to clarify aspects of compact $\sG$-manifolds of dimension $8$. Some pieces of information can be extracted from work of M.\v{C}adek et al, see \cite{Cadek-V:almostQ,Cadek-CV:quaternionic} for instance, and traced back to earlier papers \cite{Thomas:complex,Heaps:acs}, although a broader standpoint will have to await further study.

 There several ways to explain why, almost automatically,
\begin{prop}
An oriented $\sG$-manifold $\M^{8}$ is spin.
\end{prop}
\proof
Let us see three easy arguments to corroborate once more the theory's richness. First, the Stiefel-Whitney class $w_{2}(\M)\in H^{2}(\M,\Z_{2})$ is twice the Mar\-chia\-fava-Romani class \cite{Marchiafava-R:class}  for $\Sp(2)\Sp(1)$-structures, and hence null mod two.

Secondly, a quaternion-Hermitian structure renders the $8$-dimensional modules $\rS^{2}\rH\oplus\La^{3}_{0}\rE$ and $\mathrm{T}_{c}=\rE\rH$ real and $\SO(8)$-invariant; these must therefore factor through $\Sp(2)\Sp(1)$.

Thirdly, $\PSU$-structures are spin
and the irreducible  inclusion $\sG\subset\SO(5)$ lifts to $\Sp(2)=\Spin(5)\subset \Spin(8)$; in fact, the essence of Lemma \ref{psu&quat} is precisely that the fibre of this covering does not interfere with the discrete centre of $\SU(3)$.
\hfill qed\\

Instead of relying on Wu's formula to compute other Stiefel-Whitney classes in the $\Z_{2}$-cohomology ring, we shall exploit the crucial inclusion $\sG\subset \PSU$, and from  \cite{Witt:triality} we know all Stiefel-Whitney classes vanish, except possibly for $w_{4}(M)$ which squares to zero, at any rate.\\

Starting from  decomposition \eqref{splitting}
we determine the integral Pontrjagin classes $p_{i}(\M)\in H^{4i}(\M,\Z), i=1,2$.
\begin{prop}
The Pontrjagin classes of a $\sG$-manifold $\{\M^{8},g\}$  are related by
\be
\begin{split}
 \label{eq:pontrjagin}  4 p_{2}(\M)= &\; p_{1}(\M)\smile p_{1}(\M)\\
\notag p_{1}^{2}\in &\; 8640\,\Z.
\end{split}
\ee
\end{prop}
\proof
We are entitled to assume there is a circle acting on the tangent bundle  by pulling $\T_c$ back to some larger manifold fibring over $\M$, if necessary (splitting principle). Thus we can decompose
\be\label{eq:line bdls}
V_c:=V\otimes_{\R}\C=L\oplus \ol{L}\oplus\C
\ee
using a complex line bundle $L$ with Chern character $e^{x}$. Recall $\rS^{2}_{0}V_c=V_c^{\otimes 2}\ominus \La^{2}V_{c}\ominus \C\iso V_c^{\otimes 2}\ominus V_c\ominus \C$, and that $\textsl{ch}$ is a ring homomorphism from $K_{\C}(\M)$ to the even rational cohomology, to the effect that
\begin{equation*}
\begin{split}
\textsl{ch}(\T_c)& =  \textsl{ch}(V_c)+\textsl{ch}(\rS^{2}_{0}V_c)\\
& = 2\cosh x +1 + 2\cosh 2x +2\cosh x+1\\
& =  8+6x^{2}+3x^{4}.
\end{split}
\end{equation*}

\noindent
Viewing the Chern classes of the holomorphic tangent bundle $\T^{1,0}$ as elementary symmetric polynomials $s_{j}=\textstyle\sum x_{k}^{j}$ in the variables $x_{1},\ldots x_{4}$ allows to factorise the total Chern class as $ c(\T^{1,0})=\textstyle\prod(1+x_{k})$. Thus we can write the Chern character
$ch(\T^{1,0})=\textstyle\sum e^{x_{k}}=4 +s_{1}+\tfrac 12 s_{2}+\tfrac 16 s_{3}+\ldots$, whereby
  $$
  \textsl{ch}(\T_c)[\M]=\textsl{rank}\,\T_c + p_{1}+\frac 1{12}(p_{1}^{2}-2p_{2}),
  $$
   since $c_{2}=-p_{1}, c_{4}=p_{2}$, cf. \cite{Salamon:QK}.
By comparison we then have $p_{1}(\M)=6x^2, p_{2}(\M)=9x^{4}$, and the first statement
follows.

As for the second claim, it is possible to estimate the Pontrjagin number by invoking other characteristic classes. The Todd class  $\textsl{Td}(\T_c) =\textstyle\prod\frac {x_{k}}{1-e^{x_{k}}}$
equals $1+\tfrac 12 c_{1}+\tfrac 1{12}(c_{1}^{2}-p_{1})-\tfrac 1{24}c_{1}p_{1}-\tfrac 1{720}(c_{1}^{4}+4c_{1}^{2}p_{1}-\tfrac {11}4 p_{1}^{2}-c_{1}c_{3})$, telling that the square $p_{1}^{2}$ has to be a rather large integer, for it is divisible by $4\!\cdot\! 720$ at least, irrespective of the vanishing of the first Chern class, for example.
Results of \cite{Witt:triality} show that a compact oriented $\PSU$-manifold $\M^{8}$ satisfies $p_{1}^{2}\in 216\Z$; altogether, the lowest common multiple of the relevant factors is easily seen to equal $8640$. \hfill qed

\begin{remark}
Alas, the possibility that $p_{1}$ vanishes remains, and it would be valuable to know non-trivial examples with Dirac index equal zero.

The parallelisable manifold $\rS^{3}\times \rS^{5}$ (product of odd-dimensional spheres, here), despite having all classes zero, cannot be  a $\GG$-manifold: the five-sphere possesses no $\SO(3)$-structure whatsoever, neither irreducible nor standard, due to the wrong values of Kervaire's semi-characteristics \cite{Friedrich}.
\end{remark}
Still, it is interesting how the above is the best possible outcome of the Borel-Hirzebruch theory \cite{Hirzebruch:methods}:  neither resorting to Riemann-Roch, nor computing the $\hat{A}$-genus, improve this estimate. Indeed, it is known that a compact spin manifold fulfilling equation \eqref{eq:pontrjagin} has signature $\sigma=b_{4}^{+}-b_{4}^{-}$ equal $16$ times  $\hat{A}_{2}$, so the Hirzebruch-Thom signature Theorem will imply $60\sigma=p_{1}^{2}$, wherefore
 $$
 \hat{A}_{2}(\M^{8})=\tfrac{1}{16 \cdot 60}\int_{\M}p_{1}^{2}(\M^{8}).
 $$

At any rate, the existence of the quaternionic structure on $\M^{8}$ forces
$$
8e+p_{1}^{2}-4p_{2}=0,
$$
cf. \cite{Salamon:QK,Cadek-V:Sp2}, so evidently:
\begin{cor} A compact $\GG$-manifold $\{\M^{8},g\}$ has
$$
e(\T\M^{8})=0.
$$
\end{cor}
This is a useful obstruction, for it can prevent the existence of  $\sG$-structures.
\begin{exs}
The Gra{\ss}mannian $\SO(8)/\bigl(\SO(3)\times\SO(5)\bigr)$ of real, oriented three-planes in $\R^{8}$ fails the corollary, and as such it does not admit an infinitesimal $\sG$-action of our type. As we know, in fact, the denominator embeds in $\SO(8)$ in the `wrong' way.

The Wolf space ${\G_{2}}/{\SO(4)}$ is a quaternion-K\"ahler manifold of positive scalar curvature, so
its $\hat{A}$-genus is zero, whereas the Euler characteristic is not. That means it cannot carry
a $\sG$-structure, either.
\end{exs}

The corollary also falls out of the $\PSU$-side of the story, as we learn from \cite{Witt:triality}.
\\

The choice of \eqref{eq:line bdls} affects
\be\label{eq:U(1)str}
\rS^{4}\oplus \rS^{2} = (L^{2}\oplus\ol{L}^{2}\oplus L\oplus\ol{L}\oplus\C)\oplus L\oplus\ol{L}\oplus\C,
\ee
by singling out an almost complex structure. This comes from picking a point $z$ in the fibre of the twistor fibration $\mathbb{P}^{1}\hookrightarrow \mathcal{Z} \lto \M$ ($|z|=1$ reduces $\SO(3)$ to $\U(1)$, hence $\rS^{2}=\C\oplus \R$). The almost complex structure is defined by selecting a space of holomorphic tangent vector fields, and  there should be a whole $2$-parameter family thereof, depending on the choices of a line $L=L'\cos\theta+L''\sin\theta, L'\in\rS^{2}, L''\in\rS^{4}$ and of a trivial term from a similar combination of the $\C$s in \eqref{eq:U(1)str}.
By asserting that the space of $(0,1)$-forms is annihilated by
$$
\T^{1,0}=L^{2}+2\ol{L}+\C,
$$
where $L^{1/2}+L^{-1/2}=\rH$, we are fixing $J$. This gives back \eqref{eq:pontrjagin}, by the way.

There are other possible almost complex structures, two of which are fairly obvious: declaring  $L^{2}+L+L+\C$ to be
$(1,0)$-vectors defines, say, $J'$, while $L^{2}+\ol{L}+L+\C$ gives $J''$.
The latter is a quaternionic `twistor' structure, because $\rS^{3}\otimes \rH=(L^{3/2}+L^{1/2}+L^{1/2}+L^{-3/2})L^{1/2}=\T^{1,0}_{J^{''}}$.
By contrast $J'$, already met in \cite{Macia}, is non-quaternionic
as $\T^{1,0}_{J^{'}}\not=\rE\otimes L^{1/2}=L^{2}+L+\ol{L}+\C$.
There are three more basic almost complex structures obtained by complex
conjugation, ie coming from swapping $\T^{1,0}$ and $\T^{0,1}$.


\section{Relative intrinsic torsion}
\noindent
We now begin to describe the intrinsic torsion for $\sG$,
explaining in particular how the torsion spaces of the subordinate
$\G$-structures relate to each other.
Following \cite{Gray} we view two-forms as traceless, skew-symmetric matrices, then
 decompose the space  of intrinsic-torsion tensors
\[
\T\M^{8}\otimes \g^\perp\subset \T\M^{8}\otimes \Lambda^{2}\T^{*}\M^{8}
\]
 into irreducible modules under the action of $\sG$ at each point.
By equation \eqref{eq:2forms} thus, we immediately have
\begin{prop} \label{SO(3)it}
The intrinsic torsion $\tau_{\sG}$ of the $\GG$-structure is a tensor belonging in
$$
(\rS^2\oplus\rS^4)\otimes \frac{\so(8)}{\g}
    = 2\rS^{10}\oplus 5\rS^8\oplus 8\rS^6\oplus 10\rS^4\oplus 8\rS^2 \oplus3\R.
$$
This space has dimension $200$ and contains a $3$-dimensional subspace of
$\sG$-invariant tensors.\hfill qed
\end{prop}

$\GG$-invariant subspaces  will be the primary object of concern in section \ref{sec:invariant}.

\begin{definition}
For any given Lie group $\G$ containing $\sG$ we denote by $\tau^{\G}_{\sG}$
the intrinsic torsion of a $\G$-structure decomposed under $\sG$, and call it the \emph{$\G$-torsion relative to $\sG$}, or just \emph{relative $\G$-torsion},  $\sG$ being implicit most of the times.
\end{definition}
For simplicity, we disregard tensor-product signs and either juxtapose factors, or separate them
by a full stop.
The following lemma expresses the torsion spaces of
$\G=\SO(3)\times\SO(5)$-,\;$\PSU$-, and $\Sp(2)\Sp(1)$-structures in
terms of $\sG$-modules. The ensuing  Proposition \ref{prop:torsion} will show, in the same spirit of section \ref{intersection}, how  two among $\tau^{\SO(3)\times\SO(5)}_{\sG}, \tau^{\PSU}_{\sG}$ and $\tau^{\Sp(2)\Sp(1)}_{\sG}$ are enough to determine the third,
and hence the $\sG$-intrinsic torsion as well.

\begin{lemma}
Let $\G=\SO(3)\times \SO(5),  \PSU, \Sp(2)\Sp(1)$. The relative $\G$-torsion $\tau^{\G}_{\sG}$ of  $\{\M^8,g\}$ lives in the direct sum of the following modules:
\begin{center}
\begin{tabular}{l || cccccc|c}

   $\vs$                 &  $\rS^{10}$  & $\rS^8$   & $\rS^6$  & $\rS^4$ & $\rS^2$ & $\R$ & $\dim_\R$ \\[2mm]
        \hline\hline
   $\vs$     $\tau^{{\SO(3)\times\SO(5)}}_{{\sG}}$  &   1          &    3      &   5      &     6   &   5       &  2   & 120 \\[2mm]
\hline
   $\vs$     $\tau^{\PSU}_{\tiny{\sG}} $             &     2        &    4      &   6      &     8   &   6      &   2   & 158 \\[2mm]
\hline
   $\vs$     $\tau^{{\Sp(2)\Sp(1)}}_{{\sG}}$        &   1          &    3      &   5      &     6   &   5       &  2   & 120 \\[2mm]
\hline
\end{tabular}\end{center}
\end{lemma}\vspace{5mm}

\proof
$\G=\SO(3)\times\SO(5)$:\q
The six components of the $\SO(3)\times\SO(5)$-torsion are
$$
\tau_{\SO(3)\times\SO(5)}\in \Lambda^2V.W\oplus \rS^2_0V.W\oplus W \oplus V.\Lambda^2W \oplus V.\rS^2_0W\oplus V
$$
if $\T\M=V\oplus W$, see \cite{Naveira}. The terms are ordered on purpose, so that the $i$th module corresponds to Naveira's $i$th (irreducible) class $\sW_i^{AP}.$
The by-now-customary identifications $V\cong \rS^2_V, W\cong\rS^4_W$
gives the $\SO(3)\times\SO(5)$-torsion space decomposed relatively to the subgroup
$\SO(3)_V\times\SO(3)_W$
$$
 \rS^2_V\rS^4_W \oplus \rS^4_V\rS^4_W\oplus \rS^4_W\oplus \left(\rS^2_V\rS^6_W\oplus\rS^2_V\rS^4_W\right) \oplus
  \left(\rS^2_V\rS^8_W\oplus\rS^2_V\rS^4_W\right)\oplus \rS^2_V
$$
By taking $\sG=\SO(3)_V=\SO(3)_W$ and using Clebsch-Gordan we conclude.
\\

$\G=\PSU$:\q
From Theorem \ref{gig1},
 $\su(3)^\perp=2\rS^6\oplus 2\rS^2$ as  $\sG$-modules, and the claim is immediate.
\\

$\G=\Sp(2)\Sp(1)$:\q
Using the standard quaternionic notation
whereby $\Lambda^{2}_{0}\rE.\rE\ominus\rE$ is called $\mathrm{K}$, the four terms in
$$
\T_c\otimes\bigl(\ssp(2)\oplus\ssp(1)\bigr)^\perp= \rE.\rS^3\rH\oplus \mathrm{K}.\rS^3\rH\oplus  \mathrm{K}.\rH\oplus \rE.\rH
$$
correspond to the four basic classes $\sW_i^{AQH}$ of \cite{Cabrera2007}. Now identify $\rH\cong \rS^1_-$ and $\rE\cong \rS^3_+$, and decompose under $\Sp(1)_+\times\Sp(1)_-$, to the
effect that the quaternionic torsion space relative to $\Sp(1)_+\times\Sp(1)_-$ reads
$$
\rS^3_+\rS^3_- \oplus \left( \rS^7_+\rS^3_-\oplus\rS^5_+\rS^3_-\oplus\rS^1_+\rS^3_- \right)
 \oplus  \left( \rS^7_+\rS^1_-\oplus\rS^5_+\rS^1_-
\oplus\rS^1_+\rS^1_-  \right)\oplus
 \rS^3_+\rS^1_-.
$$
The final step is the diagonal identification $\Sp(1)_+=\Sp(1)_-$, that produces
the required $\sG$-modules.
\hfill qed\\

\begin{remark}
The reader interested in a description of the irreducible modules of a general $\PSU$-structure should consult \cite{Witt:triality,Puhle}, while a detailed analysis of the almost quaternion-Hermitian case can be found in  \cite{Cabrera2007}. The result about the quaternionic group was  proved in \cite{Macia}.
\end{remark}

We shall treat in the next section the two-dimensional subspace
of $\sG$-invariant tensors common to all relative torsion spaces
(cf. penultimate column in previous table).

For the last, summarising result of this part we need new labels: so let us write
$$
\rP=\SO(3)\times\SO(5),\q\rR=\PSU,\q \rQ=\Sp(2)\Sp(1),
$$
and denote with  $\mathfrak{p},\;\mathfrak{r},\;\mathfrak{q}$ the corresponding Lie
algebras. Theorem \ref{gig1} says
$$
\T^*\otimes \mathfrak{p}^\perp = \left(\T^*\otimes \frac{\mathfrak{r}}{\g}\right) \oplus
\left(\T^*\otimes\frac{\mathfrak{q}}{\g}\right).
$$
Now call   $\tau^{\rP}_\sG(\rR)\in \T^*\otimes (\mathfrak{q}/\g)$
the component of  $\tau^{\rP}_\sG$ appearing in $\tau^{\rR}_\sG$
but not present in $\tau^{\rQ}_\sG$, and similarly for
$\tau^{\rP}_\sG(\rQ)\in\T^*\otimes (\mathfrak{r}/\g)$. We can then write,
informally,
$$\label{PRQ}
\tau^{\rP}_\sG= \tau^{\rP}_\sG(\rR) \oplus \tau^{\rP}_\sG(\rQ).
$$
A similar argument makes it easy to check that these components satisfy certain
relations, for all permutations of $\rP,\rR,\rQ$, as in
\begin{prop}\label{prop:torsion}
The tensor $\tau_\sG$ of  $\{\M^8,g\}$ determines $\rP$-,\ $\rQ$-,\ $\rR$-structures
  whose relative torsion tensors $\tau^{\rP}_\sG,\tau^{\rQ}_\sG,\tau^{\rR}_\sG$
 satisfy the cyclic conditions
 \begin{eqnarray*}
\tau^{\mathrm \rP}_\sG(\rR) = \tau^{\rR}_\sG(\rP),\\
\tau^{\rP}_{\sG} = \tau^{\rP}_\sG(\rR) \oplus \tau^{\rP}_\sG(\rQ),\\
\tau_\sG = \tau^{\rP}_\sG(\rR) \oplus \tau^{\rR}_\sG(\rQ)\oplus
\tau^{\rQ}_\sG(\rP).
\end{eqnarray*}
In particular, any two yield the third.
\hfill {\textnormal{qed}}
\end{prop}


\section{Invariant torsion}\label{sec:invariant}

\noindent
It is worth remarking that $\sG$ stabilises certain exterior differential forms, as sanctioned by the
`$\R$' terms in the previous table or by the singlets in
\begin{eqnarray*}
\Lambda^3 & \cong &\mathrm{S^8\oplus 3S^6\oplus 3S^4\oplus 3S^2\oplus 2\mathbb{R}}\ \cong \ \Lambda^5,\\
\Lambda^4 & \cong & \mathrm{2S^8\oplus 2S^6\oplus 6 S^4\oplus 2S^2\oplus 2\mathbb{R}}.
\end{eqnarray*}
These $\sG$-invariant forms are two $3$-forms $\aaa,\bbb$ and one $4$-form $\gamma$,
 together with the Hodge duals $*\gamma\in\Lambda^4$,  $*\aaa,*\bbb\in\Lambda^5$.
Their algebraic nature can be described by tracking down the exact module they belong in:
 $\alpha$ appears in the decomposition of $\Lambda^3\rS^2\subset \Lambda^3,$ while
 its dual $*\alpha$ shows up in $\Lambda^5\rS^4\subset\Lambda^5.$ The form $\beta$ spans the
one-dimensional subspace in $\Lambda^2\rS^4\otimes\rS^2\subset \Lambda^3,$ and
 $*\beta$ lives in $\Lambda^2\rS^2\otimes\Lambda^3\rS^4\subset \Lambda^5.$
 Finally, $\gamma$ sits in $\Lambda^2\rS^2\otimes\Lambda^2\rS^4\subset\Lambda^4,$
whereas $*\gamma\in \rS^2\otimes \Lambda^3\rS^4\subset\Lambda^4.$

We restrict the study to $\sG$-structures with invariant intrinsic torsion, so
 the exterior differential $d:\Lambda^k\rightarrow\Lambda^{k+1}$
is a $\sG$-invariant map. The six invariant forms must then define a subcomplex of deRham's  complex,
\begin{diagram}
\La^{3}{\mathrm{T^{*}M}} & \rTo^{d} & \La^{4}{\mathrm{T^{*}M}} & \rTo^{d} &
\La^{5}{\mathrm{T^{*}M}}\\
 \uInto &   & \uInto  &   & \uInto  \\
\Span_\R\{\aaa,\bbb\} & \rTo^{d}  & \Span_\R\{\gamma,*\gamma\}  &
\rTo^{d}  & \Span_\R\{*\aaa,*\bbb\}
\end{diagram}
and the restricted $d$ on the second line is determined by linear maps between the spaces
spanned by the invariant couples. Thus the pair
 $(d\alpha,d\beta)\in 2\R\subset\Lambda^4$ must be such that
  $$
  (d\alpha,d\beta)=(\gamma,*\gamma)A,
  $$
  where $A$ is a linear transformation acting on the right on the frame $(\gamma,*\gamma)$ of $\R\oplus\R\subset\Lambda^4.$ In a completely similar manner
$$
(d\gamma,d\hodge\gamma)=(*\aaa,*\bbb)B,
$$
 $B$ being a linear map acting on $(*\aaa,*\bbb)\in 2\R\subset \Lambda^5.$

As our invariant-torsion scenario does not allow for invariant forms of  degree higher than five,
 $*\aaa$ and $*\bbb$ are forced to be closed.
 The condition $d^2=0$ becomes a non-linear constraint on the coefficients of $A,B$
\[
BA=\left[ \begin{array}{cc}
b^1_1 & b^1_2\\
b^2_1 & b^2_2\end{array}\right]
\left[ \begin{array}{cc}
a^1_1 & a^1_2\\
a^2_1 & a^2_2\end{array}\right]
=\left[\begin{array}{cc}
0 & 0\\
0& 0
\end{array}\right].
\]
In the sequel we shall extract information on the torsion out of the defining matrices. To set the pace we provide the full proof of one intermediate result only, in preparation for Theorem \ref{prop:cases}, in order to show the type of computations and arguments involved. The key point to bear in mind is that, eventually, the invariant intrinsic torsion will depend linearly on the three terms
\be\label{eq:3components}
\begin{array}{rl}
t_0 & \textnormal{from}\quad  \rS^4\otimes\rS^4 \subset \T^*\otimes (\su(3)/\g),\\[1mm]
t_+ & \textnormal{from}\quad  \tau_\sG^\rP(\rR)\in\T^*\otimes \bigl(\ssp(2)\oplus\ssp(1)\bigr)/ \g,
\\[1mm]
t_- & \textnormal{from}\quad \tau^\rQ_\sG(\rR)\in\T^*\otimes (\so(3)\oplus\so(5))/\g.
\end{array}
\ee
\begin{lemma}
The coefficients $b^1_2$ and $a^1_1$
are proportional, and contain all the information provided by the invariant torsion
in $ \tau_\sG^\rP(\rQ)\in \T^*\otimes (\su(3)/\g).$
\end{lemma}
\proof
The  differential of an arbitrary $\sG$-invariant form $\Phi$ can be
expressed by the action of the intrinsic torsion $\tau$, hence written as $d\Phi=\tau\circ\Phi$.
Let $\Pi_k(\rS^i\otimes \rS^j)$ indicate the projection to the submodule
$\rS^k$ in the Clebsch--Gordan expansion of $\rS^i\otimes\rS^j$.

We define a $\sG$-equivariant mapping $\mathrm F_{1} :
\Lambda^3\rS^2\rightarrow \Lambda^2\rS^2.\Lambda^2\rS^4$ by
skew-symmetrising the interior product of $\Phi\in\Lambda^3\rS^2$ with an element
 $X\in \tau^\rP_\sG(\rQ)=\T^*\otimes (\su(3)/\g)$
\[
   \mathrm F_{1}(\Phi) := \mathcal A(X\hook \Phi)\in \mathcal A (\rS^4\otimes
\Pi_4(\rS^4\otimes\rS^2)\otimes \rS^2\otimes\rS^2)
\]
It is easy to see this is the unique map sending $\Lambda^3\rS^2$ to $\Lambda^2\rS^2.\Lambda^2\rS^4$. Moreover, Schur's Lemma guarantees it is a
multiple of the identity map, that sends
$\alpha$ to $\gamma$.
The proportionality factor is necessarily of the form $f^\alpha t_0$, where
$f^{\alpha}\in\R\backslash \{0\}$ and $t_0$ is the  $\sG$-invariant torsion
component in $\rS^4\otimes\rS^4 \subset \T^*\otimes (\su(3)/\g)$
$$
d\alpha = \mathrm F_{1}(\alpha) = f^\alpha t_0 \gamma = a^1_1\gamma.
$$

On the same grounds there is  a unique $\sG$-equivariant function
defining $d:\rS^2.\Lambda^3\rS^4\rightarrow \Lambda^5\rS^4$, namely
$\mathrm F_{2} : \rS^2\otimes\rS^4\otimes\rS^4\otimes \rS^4 \lto \mathcal
A(\rS^4\otimes\Pi_4(\rS^4\otimes\rS^2)\otimes\rS^4\otimes\rS^4\otimes\rS^4).
$
It involves contracting with the only component of the
$\sG$-invariant intrinsic torsion in $\tau^\rP_\sG(\rQ),$ so again
 $ \mathrm F_{2}$  maps $*\gamma$ to $*\alpha$ isomorphically
$$
d(*\gamma)=f^\gamma t_0 (*\alpha),
$$
with $f^\gamma\not=0$. Eventually, $b^1_2=f^\gamma t_0$ and the proof is complete.
\hfill qed\\

We are in the position of generalising the previous discussion to all six $\sG$-invariant
forms, and define the following maps in analogy to $\mathrm F_{1}, \mathrm F_{2}$
\begin{eqnarray*}
\G_1(\rS^2\otimes\cdots)=\mathcal A (\rS^2\otimes\Pi_2(\rS^2\otimes\rS^2)\otimes\cdots),
& {} &
\rH_1(\rS^2\otimes\cdots)=\mathcal A(\rS^2\otimes\Pi_4(\rS^2\otimes\rS^2)\otimes\cdots),\\\G_2(\rS^4\otimes\cdots)=\mathcal A(\rS^2\otimes\Pi_4(\rS^2\otimes\rS^4)\otimes\cdots),
& {} &
\rH_2(\rS^2\otimes\cdots)=\mathcal A(\rS^4\otimes\Pi_2(\rS^4\otimes\rS^2)\otimes\cdots),\\
\G_3(\rS^4\otimes\cdots)=\mathcal A(\rS^4\otimes\Pi_2(\rS^4\otimes\rS^4)\otimes\cdots),
& {} &
\rH_3(\rS^4\otimes\cdots)=\mathcal A(\rS^4\otimes\Pi_4(\rS^4\otimes\rS^4)\otimes\cdots).
\end{eqnarray*}
The usual representation-theoretical argument gives
\begin{eqnarray*}
d\alpha & = &  \mathrm F_{1}(\alpha) = f^\alpha t_0 \gamma;\\[2mm]
d\beta &  = & \left(\G_1 + \G_2 + \G_3 + \rH_1+ \rH_2 + \rH_3\right)\beta \\[1mm]
       & = & \left(g^\beta_1(t_\pm)+g^\beta_2(t_\pm)+g^\beta_3t_0\right)\gamma +
\left(h^\beta_1(t_\pm) + (h^\beta_2+h^\beta_3)t_0\right)\hodge\gamma;\\[2mm]
d\gamma & = & \left(\rH_1+\rH_2 + \rH_3\right)\gamma = \left(h_1^\gamma(t_\pm) + (h^\gamma_2+h^\gamma_3)t_0\right)\hodge\beta;\\[2mm]
d(*\gamma) & = & \left( \mathrm F_{2} + \G_1 + \G_2 + \G_3\right)\hodge\gamma\\[1mm]
 &  = & \left(f^\gamma t_0\right) \hodge\alpha + \left(g^\gamma_1(t_\pm)+g^\gamma_2(t_\pm)+g^\gamma_3t_0\right)\hodge\beta,
\end{eqnarray*}
where $t_0, t_{+}, t_{-}$ are as in \eqref{eq:3components}.
As for the rest,
$f^\alpha,g^\beta_3,h^\beta_2,h^\beta_3,h^\gamma_2,h^\gamma_3,g^\gamma_3,f^\gamma$
are constants, while the remaining
$g^{\beta}_1,g_2^\beta,h_1^\beta,h_1^\gamma,g_1^\gamma,g_2^\gamma$ are
linear functions of $t_\pm$; for example
$g^\beta_1(t_\pm)=g^{\beta-}_1t_- + g^{\beta +}_1t_+.$

The entries of matrices $A, B$ depend linearly on $t_\pm,t_0$; we omit to write the explicit expressions merely for the sake of brevity.
Simple computations produce other constraints, which will not be stated formally for the same reasons. Overall, though, the picture is that the $\sG$-invariant intrinsic torsion is housed in
\begin{eqnarray*}
 A
& = &  \left[\begin{array}{c|c}
f^\alpha t_0    &   (g^{\beta -}_1+g^{\beta -}_2)t_--f^\alpha t_0 + (g_1^{\beta
  +}+g_2^{\beta +})t_+  \\
\hline
0   &    h^{\beta-}_1t_- + h_1^{\beta+}t_+
\end{array}\right]\\[3mm]
B
    & = & \left[\begin{array}{c|c}   0   &
   m f^\alpha t_0 \\
\hline
    h_1^{\gamma-}t_-+ (h_2^\gamma+h_3^\gamma)t_0 + h_1^{\gamma +}t_+  &
    (g_1^{\gamma-}+g_2^{\gamma-})t_-+g_3^\gamma t_0 -h_1^{\gamma+}t_+ \end{array}\right]
\end{eqnarray*}
subject to linear constraints
$$ \ba{c}
   f^\alpha h^{\beta\pm}_1=f^\alpha h^{\gamma\pm}_1= f^\alpha(h_2^\gamma+h_3^\gamma)=0,\quad
   (g_1^{\beta\pm}+g_2^{\beta\pm})(h_2^\gamma+h_3^\gamma)+h_1^{\beta\pm}g_3^\gamma=0,\\
   (g_1^{\beta\pm}+g_2^{\beta\pm})h_1^{\gamma\pm}+h_1^{\beta\pm}(g_1^{\gamma\pm}+g_2^{\gamma\pm})=0,\quad
   (g_1^{\beta\pm}+g_2^{\beta\pm})h_1^{\gamma\mp}+h_1^{\beta\pm}(g_1^{\gamma\mp}+g_2^{\gamma\mp})=0,
   \ea
$$ where $h_1^{\gamma+}$ is proportional to $(g_1^{\gamma +}+g_2^{\gamma+}).$
\\

To sum-up,

\begin{theorem}\label{prop:cases}
Let $\{\M^8,g\}$ be a $\sG$-manifold equipped with the six
  $\sG$-invariant forms $\{\alpha,\beta,\gamma,*\gamma,*\alpha,*\beta\}.$ If the intrinsic
  torsion is $\sG$-invariant,  the differential forms satisfy
  one of the following four sets of differential equations
\begin{center}
\begin{tabular}{l||cc|cc|cc}
     &   $d\alpha$   &  $d\beta$                     &  $d\gamma$&
  $d(*\gamma)$ &  $d(*\alpha)$ & $d(*\beta)$ \\[1mm]
\hline
$\mathrm{I}$    & $a^1_1\gamma$ & $a^1_2\gamma$                  &
$0$    &     $ m a^1_1(*\aaa)+b^2_2(*\bbb)$ & $0$ & $0$         \\[1mm]
$\mathrm{II}$    &    $ 0$         & $a^1_2\gamma + a^2_2(*\gamma)$ &
$b^2_1(*\bbb)$     &  $-((a^1_2b^2_1)/a^2_2)(*\bbb)$ & $0$ & $0$       \\[1mm]
$\mathrm{III}$    &    $ 0$         & $a^1_2\gamma$                  &
$0$   &       $b^2_2(*\bbb)$  & $0$ & $0$      \\[1mm]
$\mathrm{IV}$   &   $  0 $        &  $0$                          &
$b^2_1(*\bbb)$      &       $b^2_2(*\bbb)$  & $0$ & $0$      \\
\end{tabular}\end{center}
\end{theorem}

\proof
The aforementioned constraints on $A, B$ can be recast by the more elegant
\begin{eqnarray}
\label{short1+2}b^1_2a^2_2= m a^1_1a^2_2 = b^2_1a^1_1=0\\
\label{long}b^2_1a^1_2+b^2_2a^2_2=0
\end{eqnarray}
for some real $m$, and four generic cases ought to be considered.

I. Suppose $a^1_1\neq 0,$ so $b^1_2=ka^1_1\neq 0;$ in order
to satisfy equations \eqref{short1+2} we have to impose
$a^2_2=b^2_1=0.$ Then \eqref{long} holds necessarily,  for
any  $a^1_2,\;b^2_2,$ and the ranks $r_{A}, r_{B}$ are both equal to $1$.
At the level of forms,
\[
d\aaa  =  a^1_1\gamma,          \qquad
d\bbb  =  a^1_2 \gamma,        \qquad
d\gamma  =  0,                  \qquad
d(*\gamma)  =  m a^1_1 (*\aaa) + b^2_2 (*\bbb). \qquad
\]

II. Assume $a^1_1=0,$ so that $b^1_2=0,$ and we are left
with \eqref{long} only. Supposing $a^2_2\neq 0$ we obtain an extra
relation implying $b^2_2 =- (a^1_2b^2_1)/a^2_2$ for arbitrary
$a^1_2,\;b^2_1.$ The rank of $A$ is one, while
$r_B=0 \textnormal{\ or\ } 1$ depending on $b^2_1.$ Therefore
\[
d\aaa  =  0,\qquad
d\bbb  =  a^1_2\gamma + a^2_2(*\gamma),\qquad
d\gamma  =  b^2_1 (*\bbb),\qquad
d(*\gamma) =  -\frac{a^1_2b^2_1}{a^2_2} (*\bbb).
\]

III. For $a^1_1=0$ and $a^2_2=0,$ we are left with $b^2_1a^1_2=0.$
Taking $a^1_2\neq 0$ forces $b^2_1=0$, and
$b^2_2$ is free. Again, $r_{A}=1$, and $b^2_2$ decides whether $r_B=0 \textnormal{\ or\ } 1$.
The forms satisfy
\[
d\aaa  =   0,\qquad
d\bbb  =  a^1_2 \gamma,\qquad
d\gamma =  0,\qquad
d(*\gamma) =  b^2_2 (*\bbb).
\]

IV. If in the previous case we assume $a^1_2=0$ then $b^2_1$
becomes a free parameter, as does $b^2_2.$ Now
$A$ is null and $0\leq r_B \leq 1$ depending on $b^2_1,b^2_2.$ Hence\smallbreak
\[
d\aaa  =  0,\qquad
d\bbb  =  0,\qquad
d\gamma  =  b^2_1(*\bbb),\qquad
d(*\gamma) =  b^2_2 (*\bbb).
\]
\hfill  qed
\vspace{1cm}

One fact deserves an explanation.
Since the quaternionic $4$-form $\Omega$  of Kraines and Bonan \cite{Kraines:topology,Bonan:algebre-exterieure} can be decomposed  as
$$
\Omega=\gamma+\hodge\gamma
$$
using a suitable four-form $\gamma$, and its derivative
$d(\gamma + *\gamma) = b^1_2(*\alpha) +(b^2_1+b^2_2)(*\beta)$
 `is' the quaternionic torsion $\tau_{\Sp(2)\Sp(1)}$ in disguise, one can
 easily see that the coefficients
are linear in $\{t_-,t_0\}$ only, as $t_+$  does not appear in the
invariant relative $\Sp(2)\Sp(1)$-torsion.

On the other hand there exists a second $4$-form
\be\label{Omega'}
\Omega'=\gamma-\hodge\gamma
\ee
with  stabiliser  $\Sp(2)\Sp(1)\subset\GL(2,\H)\H^{*}$, that arises
by changing the orientation of the almost quaternion-Hermitian structure. Unlike
 $\Omega$, the anti-selfdual $\Omega'$ is a function of $t_0$ and of a linear combination of both $t_+, t_{-}$. The presence of two (non-conjugated) Lie groups isomorphic to $\Sp(2)\Sp(1)$ is consistent with Remark \ref{remtri1}.
\\

\begin{exs}
To finish, we indicate how to find examples befitting the Theorem's cases.

\textnormal{I}:
Both $A$ and $B$ have rank $1$, so we can
take $\M=\SU(3)$ as in example \ref{ex:oscar}.
\smallbreak

\textnormal{II}:
While $A$ has rank $1$, for $B$ it is either $0$ or $1$. Since
the quaternionic form $\Omega$ is closed iff $b^2_1+b^2_2=0$,
necessarily
  $\Omega=d\beta$ is exact, and many constructions are known, see
 \cite{Swann, Salamon:Bilbao}.\smallbreak

\textnormal{III}:
 Suppose we insist on wanting $\Omega=\gamma + *\gamma$ closed. As $\gamma$ is always exact, we must force
 $B$ to be null and therefore there are no $\GG$-invariant  four-forms.
\smallbreak

\textnormal{IV}: Invariant three-forms are closed ($A=0$), and \cite{GM,Miquel} provide us with  constructed through special foliations.

Further examples are easy to build, and the previous computations indicate that a classification of sorts is within sight. The $5$-dimensional theory has the advantage of highlighting the prominent geometrical aspects of dimension $8$, and one can take, for instance, a rank-three bundle over one explicit $5$-dimensional Lie group of \cite{FinoChiossi}, or use Cartan-K\"ahler techniques as in the last section of \cite{Nurowski}.  The instances of \cite{Friedrich2} suggest nice twistor-flavoured constructions of similar type.
\end{exs}

\end{document}